\documentclass[3p,12pt]{elsarticle}

\usepackage{fleqn}

\usepackage{geometry}
\usepackage{natbib}
\usepackage{pifont}
\usepackage{pifont}
\usepackage{natbib}
\usepackage{amssymb}

\usepackage{amsmath}
\newtheorem{theorem}{Theorem}[section]
\newtheorem{lemma}[theorem]{Lemma}

\newtheorem{algorithm}[theorem]{Algorithm}
\newtheorem{remark}[theorem]{Remark}

\numberwithin{equation}{section}
\journal{Journal of Nonlinear and Convex Analysis}
\begin{document}

\begin{frontmatter}

\title{ Metric projection and convergence theorems for  nonexpansive mappings in Hadamard spaces}
%A characterization of metric projection in  Hadamard spaces with applications }
%Strong convergence of iterative methods for nonexpansive mappings in Hadamard spaces by metric projection }

%% use optional labels to link authors explicitly to addresses:
%% \author[label1,label2]{<author name>}
%% \address[label1]{<address>}
%% \address[label2]{<address>}

\author[]{Hossein Dehghan\corref{cor1}}
\cortext[cor1]{Corresponding author.}
\ead{h$_{-}$dehghan@iasbs.ac.ir, hossein.dehgan@gmail.com}

\author{Jamal Rooin}
\ead{rooin@iasbs.ac.ir}
\address{Department of Mathematics, Institute for Advanced Studies in Basic Sciences (IASBS), Gava Zang, Zanjan 45137-66731, Iran
 }

%\cortext[cor1]{Corresponding author.}

\begin{abstract}
%% Text of abstract
For a nonempty convex subset $C$ of a Hadamard space $X$, it is proved that $u=P_Cx$ if and only if $\langle\overrightarrow{xu},\overrightarrow{uy}\rangle \geqslant0$ for all $y\in C$. As an application of this characterization, we prove strong convergence of two iterative algorithms with perturbations for
nonexpansive mappings.
\end{abstract}

\begin{keyword}
%% keywords here, in the form: keyword \sep keyword
% MSC codes here, in the form: \MSC code \sep code
Hadamard space; metric projection; quasilinearization; nonexpansive mapping; iterative method.\\
 \MSC[2010] 47H09 \sep 47H10 % 47J20 Variational and other types of inequalities involving nonlinear operators
\end{keyword}

\end{frontmatter}

%%
%% Start line numbering here if you want
%%
% \linenumbers

%% main text

\section{Introduction }

\par
A metric space $(X,d)$ is a CAT(0) space if it is geodesically connected and if every geodesic triangle in $X$ is at least as thin as its
comparison triangle in the Euclidean plane. For other equivalent definitions and basic properties, we refer the reader to standard texts such as \cite{Ballmann,Bridson}.
Complete CAT(0) spaces are often called Hadamard spaces.
 % (see \cite{Kirk R-tree}).
Let $x,y\in X$ and $\lambda\in [0,1]$. We write $\lambda x\oplus (1-\lambda) y$ for the unique point $z$ in the geodesic segment joining from $x$ to $y$ such that
\begin{eqnarray}\label{oplus}
d(z, x)= (1-\lambda)d(x, y)\ \ \ \mbox{and}\ \ \  d(z, y)=\lambda d(x, y).
\end{eqnarray}
We also denote by $[x, y]$ the geodesic segment joining from $x$ to $y$, that is, $[x, y]=\{\lambda x\oplus (1-\lambda) y : \lambda\in [0, 1]\}$. A subset $C$ of a CAT(0) space is convex if $[x, y]\subseteq C$ for all $x, y\in C$.
\par
%It is well known that a normed linear space satisfies CAT(0)-inequality if and only if it is a pre-Hilbert space, hence it is not so unusual to have an inner product-like notion in Hadamard spaces.
Berg and Nikolaev in \cite{Berg1} introduced the concept of \emph{quasilinearization} in a metric space $X$. Let us formally denote a pair $(a,b)\in X\times X$ by $\overrightarrow{ab}$ and call it a vector. Then quasilinearization is a map $\langle \cdot, \cdot \rangle : (X\times X)\times(X\times X)\to \mathbb{R}$ defined by
\begin{eqnarray}\label{qasilin}
 \langle\overrightarrow{ab},\overrightarrow{cd}\rangle = \frac{1}{2}\left(d^2(a,d)+d^2(b,c)-d^2(a,c)-d^2(b,d)\right), \ \ \ \ (a,b,c,d\in X).
\end{eqnarray}
It is easily seen that $\langle\overrightarrow{ab},\overrightarrow{cd}\rangle =\langle\overrightarrow{cd},\overrightarrow{ab}\rangle$,  $\langle\overrightarrow{ab},\overrightarrow{cd}\rangle = -\langle\overrightarrow{ba},\overrightarrow{cd}\rangle$ and $\langle\overrightarrow{ax},\overrightarrow{cd}\rangle + \langle\overrightarrow{xb},\overrightarrow{cd}\rangle = \langle\overrightarrow{ab},\overrightarrow{cd}\rangle$ for all $a,b,c,d, x\in X$.
We say that $X$ satisfies the Cauchy-Schwarz inequality if
\begin{eqnarray}%\label{Cha-Sch}
\nonumber \langle\overrightarrow{ab},\overrightarrow{cd}\rangle\leqslant d(a,b) d(c,d)
\end{eqnarray}
for all $a,b,c,d\in X$. It is known \cite[Corollary 3]{Berg1} that a geodesically connected metric space is a CAT(0) space if and only if it satisfies the Cauchy-Schwarz inequality.

%%%%%%%%%%%%%%%%%%%%%%%%%%%%%%%%%%%%%%%%%%%%%%%%%%%%%%%%%%%%%%%%%%%%%%%%%%%%%%%%%%%%%%%%%%%%%%%%%%%%%%%%%%%%%%%%%%%%%%%%%%%%%%%%%%%%%%%%%%%%%%%%%%%%%%%%

\par
Let $C$ be a nonempty subset of a complete CAT(0) space $X$. Then a mapping $T$ of $C$
into itself is called nonexpansive if $d(Tx, Ty)\leqslant d(x, y)$ for all $x, y\in C$. A point $x \in C$ is called a fixed point of $T$ if $Tx=x$. We denote by $F(T)$ the set of all fixed points of $T$. Kirk  \cite{Kirk1} showed that the fixed point set of a nonexpansive mapping $T$ is nonempty, closed and convex.\\
Iterative methods for finding fixed points of nonexpansive mappings have
received vast investigations due to its extensive applications in a variety of applied
areas of inverse problem, partial differential equations, image recovery, and signal processing;
see \cite{Xu,Wu-Chang-Yuan,Marino-Xu,Shimoji-Takahashi,Petrusel-Yao,Chidume-Chidume,Yao-Shahzad} and the references therein. One of the difficulties in carrying out results from Banach space to Hadamard space setting lies in the heavy use of the linear structure of the Banach spaces.
% Yao and Shahzad \cite{Yao-Shahzad} suggested new iterative methods with perturbations in Hilbert space and proved convergence theorems of the proposed iterative algorithms.
\par
Now having an inner product-like notion( quasilinearization) in Hadamard spaces,
%it is natural to ask: Do the results of Yao and Shahzad still hold in Hadamard spaces? \par In this paper, the above concern is answered in the affirmative in Hadamard space setting.
we first obtain a characterization of
metric projection together with some basic lemmas in Hadamard
spaces. Then, we introduce two iterative methods to approximate
fixed points of nonexpansive mappings in Hadamard spaces.

\section{Preliminaries and lemmas}

In this section, we recall some basic results and prove some useful lemmas which we need in the sequel.
%To proceed in this direction, we need the following lemmas.

\begin{lemma}\label{hyberbolic ineqal}\cite[Proposition 2.2]{Bridson}
Let $X$ be a CAT(0) space, $p, q, r, s\in X$ and $\lambda\in [0,1]$. Then
\begin{eqnarray}
\nonumber d(\lambda p\oplus (1-\lambda) q, \lambda r\oplus (1-\lambda) s)\leqslant \lambda d(p,r)+(1-\lambda)d(q,s).
\end{eqnarray}
\end{lemma}

%%%%%%%%%%%%%%%%%%%%%%%%%%%%%%%%%%%%%%%%%%%%%%%%%%%%%%%%%%%%%%%%%%%%%%%%%%%%%%%%%%%%%%%%%%%%%%%%%%%%%%%%%%%%%%

\begin{lemma}\label{triangle ine}\cite[Lemma 2.4]{Dhompongsa}
Let $X$ be a CAT(0) space, $x,y,z\in X$ and $\lambda\in [0,1]$. Then
\begin{eqnarray}
\nonumber d(\lambda x\oplus (1-\lambda) y, z)\leqslant \lambda d(x,z)+(1-\lambda)d(y,z).
\end{eqnarray}
\end{lemma}

\begin{lemma}\label{para ine}\cite[Lemma 2.5]{Dhompongsa}
Let $X$ be a CAT(0) space, $x,y,z\in X$ and $\lambda\in [0,1]$. Then
\begin{eqnarray}
\nonumber d^2(\lambda x\oplus (1-\lambda) y, z)\leqslant \lambda d^2(x,z)+(1-\lambda)d^2(y,z)-\lambda(1-\lambda)d^2(x,y).
\end{eqnarray}
\end{lemma}

%%%%%%%%%%%%%%%%%%%%%%%%%%%%%%%%%%%%%%%%%%%%%%%%%%%%%%%%%%%%%%%%%%%%%%%%%%%%%%%%%%%%%%%%%%%%%%%%%%%%%%%%%%%%%%%%%%%%%%%%%%%%%%%%%%%%%%%%%%%%%%%%%%%%%%%%
 \begin{lemma}\cite[Lemma 1.1]{Panyanak-Cuntavepanit}\label{suzuki lem}
Let $\{x_n\}$, $\{y_n\}$ and $\{z_n\}$ be bounded sequences in a metric space of hyperbolic type $X$ and $\{\beta_n\}$ be a sequence in $[0, 1]$ with $0 < \liminf_{ n\to\infty} \beta_n\leqslant \limsup_{ n\to\infty} \beta_n  < 1$. Suppose that $x_{n+1} =  (1 - \beta_n)y_n \oplus \beta_nz_n $ for all $n\geqslant1$, $\lim_{ n\to\infty} d(y_n, x_n) = 0$ and $\limsup_{ n\to\infty} ( d(z_{n+1}, z_n) - d(x_{n+1}, x_n)) \leqslant0$. Then, $\lim_{ n\to\infty} d(z_n, x_n) = 0$.
\end{lemma}
We note that every CAT(0) space is of hyperbolic type (see \cite{Kirk}).
%The proof closely follows the proof of Lemma 2.1 and 2.2 in \cite{Suzuki}. Therefore, we omit the details.
%%%%%%%%%%%%%%%%%%%%%%%%%%%%%%%%%%%%%%%%%%%%%%%%%%%%%%%%%%%%%%%%%%%%%%%%%%%%%%%%%%%%%%%%%%%%%%%%%%%%%%%%%%%%%%
\begin{lemma}\label{Liu' lem}
(Liu's lemma) Assume that $\{a_n\}$ is a sequence of nonnegative real numbers such that
\begin{eqnarray}
\nonumber a_{n+1}\leqslant (1-\gamma_n) a_n + \gamma_n\delta_n + \sigma _n,\ \ \ \ n\geqslant 0,
\end{eqnarray}
where $\{\gamma_n\}$ is a sequence in (0,1), $\{\delta_n\}$ is sequence in $\mathbb{R}$ and $\{\sigma_n\}$ is a sequence of nonnegative numbers such that (i) $\lim_{n\to \infty} \gamma_n=0$ and $\sum_{n=0}^\infty \gamma_n=\infty$, (ii) $ \limsup_{n\to \infty} \delta_n \leqslant 0$ or $\sum_{n=0}^\infty \gamma_n|\delta_n|<\infty$, (iii) $\sum_{n=0}^\infty \sigma_n<\infty$. Then $\lim_{n\to \infty}a_n=0$.
\end{lemma}

We shall repeatedly use the following useful lemmas in the next sections.
%%%%%%%%%%%%%%%%%%%%%%%%%%%%%%%%%%%%%%%%%%%%%%%%%%%%%%%%%%%%%%%%%%%%%%%%%%%%%%%%%%%%%%%%%%%%%%%%%%%%%%%%%%%%%%%%%%%%%%%%%%%%%%%%%%%%%%%%%%%%%%
\begin{lemma}\label{coffi quasi lem}
Let $X$ be a CAT(0) space, $x,y\in X$, $\lambda\in [0,1]$ and $z= \lambda x\oplus (1-\lambda)y$. Then,
\begin{eqnarray}
\nonumber \langle\overrightarrow{zy},\overrightarrow{zw}\rangle \leqslant \lambda \langle\overrightarrow{xy},\overrightarrow{zw}\rangle\ \ \ (w\in X).
\end{eqnarray}
\end{lemma}
\textbf{Proof.} Using (\ref{oplus}) and Lemma \ref{para ine}, we have
\begin{eqnarray}
\nonumber 2(\langle\overrightarrow{zy},\overrightarrow{zw}\rangle-\lambda \langle\overrightarrow{xy},\overrightarrow{zw}\rangle)&=& d^2(z,w) + d^2(y,z) - d^2(y,w)\\
\nonumber && -\lambda \left( d^2(x,w) + d^2(y,z) - d^2(x,z) - d^2(y,w)  \right)\\
\nonumber&\leqslant& \lambda d^2(x,w) + (1-\lambda) d^2(y,w) -\lambda(1- \lambda) d^2(x,y)+  d^2(y,z) \\
\nonumber && - d^2(y,w) -\lambda \left( d^2(x,w) + d^2(y,z) - d^2(x,z) - d^2(y,w)  \right)\\
\nonumber &=& (1-\lambda) d^2(y,z) + \lambda d^2(x,z) -\lambda(1- \lambda) d^2(x,y)\\
\nonumber &=& (1-\lambda) \lambda^2 d^2(y,x) + \lambda (1-\lambda)^2 d^2(x,y) -\lambda(1- \lambda) d^2(x,y)\\
\nonumber &=& 0,
\end{eqnarray}
which implies the desired inequality. $\ \ \ \Box$\\
%%%%%%%%%%%%%%%%%%%%%%%%%%%%%%%%%%%%%%%%%%%%%%%%%%%%%%%%%%%%%%%%%%%%%%%%%%%%%%%%%%%%%%%%%%%%%%%%%%%%%%%%%%%%%%%%%%%%%%%%%%%%%%%%%%%%%%%%%%%%%%%%%%%
\begin{lemma}\label{coffi quasi lem 2}
Let $X$ be a CAT(0) space, $x,y, z\in X$ and $\lambda\in [0,1]$. Then,
\begin{eqnarray}%\label{coffi quasi 2}
\nonumber d^2(\lambda x\oplus (1-\lambda)y, z)\leqslant \lambda^2 d^2(x, z) + (1-\lambda)^2 d^2(y, z) + 2\lambda(1-\lambda)\langle\overrightarrow{xz},\overrightarrow{yz}\rangle.
\end{eqnarray}
\end{lemma}
\textbf{Proof.} Using Lemma \ref{para ine}, we have
\begin{eqnarray}
\nonumber d^2(\lambda x\oplus (1-\lambda)y, z)&\leqslant& \lambda d^2(x, z) + (1-\lambda) d^2(y, z)- \lambda(1-\lambda) d^2(x,y)\\
\nonumber &=& \lambda^2 d^2(x, z) + (1-\lambda)^2 d^2(y, z) \\
\nonumber &&+ \lambda(1-\lambda)\left(d^2(x,z)+ d^2(y,z)-d^2(x,y)\right)\\
\nonumber &=& \lambda^2 d^2(x, z) + (1-\lambda)^2 d^2(y, z)+ 2 \lambda(1-\lambda) \langle\overrightarrow{xz},\overrightarrow{yz}\rangle,
\end{eqnarray}
which is the desired inequality. $\ \ \ \Box$\\
\par
Let $\{x_n\}$ be a bounded sequence in a complete CAT(0) space $X$. For $x\in X$, we set
\begin{eqnarray}
\nonumber r(x, \{x_n\})=\limsup_{n\to\infty}d(x, x_n).
\end{eqnarray}
The asymptotic radius $r(\{x_n\})$ of $\{x_n\}$ is given by
\begin{eqnarray}
\nonumber r(\{x_n\})=\inf\{r(x, \{x_n\}): x\in X\},
\end{eqnarray}
and the asymptotic center $A(\{x_n\})$ of $\{x_n\}$ is the set
\begin{eqnarray}
\nonumber A(\{x_n\})=\{x\in X : r(x, \{x_n\})=r(\{x_n\})\}.
\end{eqnarray}

It is known from Proposition 7 of \cite{Dhompongsa-Kirk-Sims} that in a CAT(0) space, $A(\{x_n\})$ consists of exactly one point.\\
A sequence $\{x_n\}\subset X$ is said to $\Delta$-converge to $x\in X$ if $A(\{x_{n_k}\})=\{x\}$ for every subsequence $\{x_{n_k}\}$ of $\{x_n\}$.
%Uniqueness of asymptotic center implies that CAT(0) space $X$ satisfies Opial's property, i.e., for given $\{x_n\}\subset X$ such that $\{x_n\}$ $\Delta$-converges to $x$ and given $y\in X$ with $y\not =x$,
%\begin{eqnarray}
%\nonumber \limsup_{n\to\infty}d(x_{n},x) < \limsup_{n\to\infty}d(x_{n}, y).
%\end{eqnarray} Since it is not possible to formulate the concept of demiclosedness in a CAT(0) setting, as stated in linear spaces, let us formally say that "$I-T$ is demiclosed at zero" if the conditions, $\{x_n\}\subseteq C$ $\Delta$- converges to $x$ and $d(x_n, Tx_n)\to 0$ imply $x\in F(T)$.
%%%%%%%%%%%%%%%%%%%%%%%%%%%%%%%%%%%%%%%%%%%%%%%%%%%%%%%%%%%%%%%%%%%%%%%%%%%%%%%%%%%%%%%%%%%%%%%%%%%%%%%%%%%%%%%%%%%%%%%%%%%%%%%%%%%%%%%%%%%%%%%%%%%%%%%%
\par
We need the following lemmas.
\begin{lemma}\label{delta convergent subseq}
{\rm \cite{Kirk-Panyanak}} Every bounded sequence in a complete CAT(0) space always has a $\Delta$-convergent subsequence.
\end{lemma}
%%%%%%%%%%%%%%%%%%%%%%%%%%%%%%%%%%%%%%%%%%%%%%%%%%%%%%%%%%%%%%%%%%%%%%%%%%%%%%%%%%%%%%%%%%%%%%%%%%%%%%%%%%%%%%%%%%%%%%%%%%%%%%%%%%%%%%%%%%%%%%%%%%%%%%%%
\begin{lemma}\label{deta limit is in C}
{\rm \cite{Dhompongsa-Kirk-Panyanak}} If $C$ is a closed convex subset of a complete CAT(0) space and if $\{x_n\}$ is a bounded sequence in $C$, then the asymptotic center
of $\{x_n\}$ is in $C$.
\end{lemma}

\begin{lemma}\label{demiclosed of T}
{\rm \cite{Dhompongsa-Kirk-Panyanak}} If $C$ is a closed convex subset of $X$ and $T:C\to X$ is a nonexpansive mapping, then the conditions $\{x_n\}$ $\Delta$-convergence to $x$ and $d(x_n,Tx_n)\to 0$, and imply $x\in C$ and $Tx=x$.
\end{lemma}
%A metric space $(X,d)$ is a CAT(0) space if it is geodesically connected and if every geodesic triangle in $X$ is at least as thin as its
%comparison triangle in the Euclidean plane. For other equivalent definitions and basic properties, we refer the reader to standard texts such as \cite{Ballmann,Bridson, Burago}.
%%Complete CAT(0) spaces are often called Hadamard spaces.
% % (see \cite{Kirk R-tree}).
%Let $x,y\in X$ and $\lambda\in [0,1]$. We write $\lambda x\oplus (1-\lambda) y$ for the unique point $z$ in the geodesic segment joining from $x$ to $y$ such that
%\begin{eqnarray}\label{oplus}
%d(z, x)= (1-\lambda)d(x, y)\ \ \ \mbox{and}\ \ \  d(z, y)=\lambda d(x, y).
%\end{eqnarray}
%We also denote by $[x, y]$ the geodesic segment joining from $x$ to $y$, that is, $[x, y]=\{\lambda x\oplus (1-\lambda) y : \lambda\in [0, 1]\}$. A subset $C$ of a CAT(0) space is convex if $[x, y]\subseteq C$ for all $x, y\in C$.
\par
%It is well known that a normed linear space satisfies CAT(0)-inequality if and only if it is a pre-Hilbert space, hence it is not so unusual to have an inner product-like notion in Hadamard spaces.
\begin{lemma}\label{delta iff limsup}{\rm \cite[Theorem 2.6]{Kakavandi}} Let $X$ be a complete CAT(0) space, $\{x_n\}$ be a sequence in $X$ and $x \in X$. Then $\{x_n\}$ $\Delta$-converges to $x$ if and only if $\limsup_{n\to\infty}\langle \overrightarrow{x x_n},\overrightarrow{x y} \rangle \leq 0$  for all $y\in X$.
\end{lemma}

%%%%%%%%%%%%%%%%%%%%%%%%%%%%%%%%%%%%%%%%%%%%%%%%%%%%%%%%%%%%%%%%%%%%%%%%%%%%%%%%%%%%%%%%%%%%%%%%%%%%%%%%%%%%%%%%%
\section{Metric projection }
\par
Let $C$ be a nonempty complete convex subset of a CAT(0) space $X$. It is known that for any $x\in X$ there exists a unique point
$u\in C$ such that
\begin{eqnarray}
\nonumber d(x,u)= \min_{ y\in C} d(x, y).
\end{eqnarray}
The mapping $P_C : X \rightarrow C$ defined by $P_Cx= u$ is called the \emph{metric projection} from $X$ onto $C$. Also, $P_C$ is nonexpansive ( see \cite[Proposition 2.4]{Bridson}).
%Recently, the authors presented the following characterization of metric projection in extended abstract of a conference \cite{Dehghan-Rooin}. For convenience of the readers, we include the complete proof.
Now, we state and prove our characterization of metric projection.
%%%%%%%%%%%%%%%%%%%%%%%%%%%%%%%%%%%%%%%%%%%%%%%%%%%%%%%%%%%%%%%%%%%%%%%%%%%%%%%%%%%%%%%%%%%%%%%%%%%%%%%%%%%%%%%%%%%%%%%%%%%%%%%%%%%%%%%5
%%%%%%%%%%%%%%%%%%%%%%%%%%%%%%%%%%%%%%%%%%%%%%%%%%%%%%%%%%%%%%%%%%%%%%%%%%%%%%%%%%%%%%%%%%%%%%%%%%%%%%%%%%%%%%%%%%%%%%%%%%%%%%%%%%%%%%%5
\begin{theorem}\label{proj}
Let $C$ be a nonempty convex subset of a CAT(0) space $X$, $x\in X$ and $u\in C$. Then  $u=P_Cx$ if and only if
\begin{eqnarray}%\label{charac}
\nonumber \langle\overrightarrow{xu},\overrightarrow{uy}\rangle\geqslant0\ \ \ \ (\forall y\in C).
\end{eqnarray}
\end{theorem}
\noindent\textbf{Proof.} Let $\langle\overrightarrow{xu},\overrightarrow{uy}\rangle\geqslant0$  for all $y\in C$. If $d(x,u)=0$, then the assertion is clear. Otherwise, we have
\begin{eqnarray}
\nonumber
\langle\overrightarrow{xu},\overrightarrow{xy}\rangle -\langle\overrightarrow{xu},\overrightarrow{xu}\rangle= \langle\overrightarrow{xu},\overrightarrow{uy}\rangle\geqslant0.
\end{eqnarray}
This together with the Cauchy-Schwarz inequality implies that
\begin{eqnarray}
\nonumber d^2(x,u)=\langle\overrightarrow{xu},\overrightarrow{xu}\rangle \leqslant \langle\overrightarrow{xu},\overrightarrow{xy}\rangle
\leqslant  d(x,u) d(x,y).
\end{eqnarray}
That is, $d(x,u)\leqslant d(x,y)$ for all $y\in C$ and so $u=P_Cx$.
\par
For the converse, let $u=P_Cx$. Since $C$ is convex, then $z=\lambda y\oplus(1-\lambda)u\in C$ for all $y\in C$ and $\lambda\in(0,1)$. Thus, $d(x,u)\leqslant d(x,z)$. Using (\ref{qasilin}) we have
\begin{eqnarray}\label{char ineq1}
 \langle\overrightarrow{xz},\overrightarrow{uz}\rangle \geqslant \frac{1}{2}d^2(x,z)-\frac{1}{2}d^2(x,u)\geqslant 0.
\end{eqnarray}
On the other hand, by using Lemma \ref{coffi quasi lem}, we have $\langle\overrightarrow{xz},\overrightarrow{uz}\rangle\leqslant \lambda\langle\overrightarrow{xz},\overrightarrow{uy}\rangle$. This together with (\ref{char ineq1}) implies that $\langle\overrightarrow{xz},\overrightarrow{uy}\rangle \geqslant 0$. Since the function $d(\cdot , x): X \to \mathbb{R}$ is continuous for all $x\in X$, considering (\ref{qasilin}) and letting $\lambda\to 0^+$, we have $\langle\overrightarrow{xu},\overrightarrow{uy}\rangle\geqslant0$. This completes the proof.
$\Box$

%%%%%%%%%%%%%%%%%%%%%%%%%%%%%%%%%%%%%%%%%%%%%%%%%%%%%%%%%%%%%%%%%%%%%%%%%%%%%%%%%%%%%%%%%%%%%%%%%%%%%%%%%%%%%%%%%%%%%%%%%%%%%%%%%%%%%%%%%%%%%%%%%%%%%%%%

\section{Convergence theorems}
In this section, we apply the obtained results to approximate fixed points of nonexpansive mappings in Hadamard spaces. In the rest of the paper, $(X, d)$ is a Hadamard space,  $o$ is an arbitrary fixed element in $X$, which we may call the "zero" of $X$ and $\|x\|:=d (x,o)$ for all $x\in X$.

\subsection{Convergence of an implicit algorithm\\}

Let $T$ be a nonexpansive self-mapping of a nonempty closed convex subset $C$ of a Hadamard space $X$. We denote by $F(T)$ the set of all fixed points of $T$. Fix $u\in X$. Then for each $\alpha\in (0, 1)$, there exists a unique point $x_\alpha\in C$ satisfying $x_\alpha= P_C\left(\alpha u \oplus (1-\alpha)Tx_\alpha\right)$ because the mapping $x\mapsto P_C\left(\alpha u \oplus (1-\alpha)Tx\right)$ is contractive by virtue of Lemma \ref{hyberbolic ineqal} and nonexpansiveness of $P_C$. Therefore, we may define the following implicit iterative method.
\begin{algorithm}
 Let $\{\alpha_m\}$ be a sequence in $(0, 1)$ and define the iterative sequence $\{x_m\}$ by
\begin{eqnarray}\label{mann like}
x_m=P_C\left(\alpha_m u_{m}\oplus(1-\alpha_m)Tx_m\right),\ \ \ \ m\geqslant1,
\end{eqnarray}
where the sequence $\{u_m\}\subset X$ is a small perturbation for the $m$-step iteration satisfying $\|u_m\|\to 0$ as $m\to \infty$.
\end{algorithm}
%\par In order to prove convergence of $\{x_n\}$ we need the following lemmas.
%%%%%%%%%%%%%%%%%%%%%%%%%%%%%%%%%%%%%%%%%%%%%%%%%%%%%%%%%%%%%%%%%%%%%%%%%%%%%%%%%%%%%%%%%%%%%%%%%%%%%%%%%%%%%%%%%%%%%%%%%%%%%%%%%%%%%%%%%%%%%%%%%%%%%%%%
\begin{theorem}
If $F(T)\not=\emptyset$, then as $\alpha_m\to 0$, the sequence $\{x_m\}$ generated by the implicit method ( \ref{mann like}) converges to a $q\in F(T)$.
\end{theorem}
\textbf{Proof.}
We first show that $\{x_m\}$ is bounded. Taking $p\in F(T)$ and using the fact that $P_C$ is nonexpansive and Lemma \ref{triangle ine}, we have
\begin{eqnarray}
\nonumber d(x_m,p)&=&d\left(P_C\left(\alpha_m u_{m}\oplus(1-\alpha_m)Tx_m\right),p\right)\\
\nonumber &\leqslant&d\left(\alpha_m u_{m}\oplus(1-\alpha_m)Tx_m,p\right)\\
\nonumber &\leqslant&\alpha_m d(u_{m}, p) +(1-\alpha_m)d(Tx_m,p)\\
\nonumber &\leqslant&\alpha_m d(u_{m}, p) +(1-\alpha_m)d(x_m,p),
\end{eqnarray}
which implies that
\begin{eqnarray}
\nonumber d(x_m,p)\leqslant d(u_{m}, p)\leqslant \|u_m\|+\|p\|.
\end{eqnarray}
Since $\|u_m\|\to 0$, then $\{u_m\}$ is bounded. It follows that $\{x_m\}$ is bounded, so is the sequence $\{Tx_m\}$. Thus, there exists a constant $M > 0$ such that $\|u_m\|, \|x_m\|, \|Tx_m\|\leqslant M$ for all $m\geqslant1$. Since $Tx_m\in C$, we get
\begin{eqnarray}\label{d(x_m, Tx_m) to 0}
 \nonumber d(x_m,Tx_m)&=&d\left(P_C\left(\alpha_m u_{m}\oplus(1-\alpha_m)Tx_m\right), P_C(Tx_m)\right)\\
\nonumber &\leqslant& d\left(\alpha_m u_{m}\oplus(1-\alpha_m)Tx_m, Tx_m\right)\\
\nonumber  &=&\alpha_m d(u_{m},Tx_m)\\
 &\leqslant&2\alpha_m M \to 0,
\end{eqnarray}
as $m\to \infty$. Setting $y_m=\alpha_m u_{m}\oplus(1-\alpha_m)Tx_m$ for all $m\geqslant1$, we then have $x_m=P_Cy_m$. Also, by Lemma \ref{triangle ine}
\begin{eqnarray}\label{d(x_m, y_m) to 0}
 d(x_m,y_m) \leqslant \alpha_m d (x_m, u_{m}) + (1-\alpha_m) d (x_m, Tx_m) \to 0,\ \ \ (as\ \ m\to\infty).
\end{eqnarray}
%It follows from (\ref{charac}) that \begin{eqnarray} \nonumber\langle\overrightarrow{x_my_m},\overrightarrow{x_mp}\rangle\leqslant 0. \end{eqnarray}
Since $P_C$ is nonexpansive and $p,Tx_m\in C$, we have
\begin{eqnarray}\label{x_m to y_m}
\nonumber 2\langle\overrightarrow{x_mTx_m},\overrightarrow{x_mp}\rangle &=& d^2(x_m,p)+ d^2(Tx_m,x_m)- d^2(Tx_m,p)\\
\nonumber &\leqslant& d^2(y_m,p)+ d^2(Tx_m,y_m)- d^2(Tx_m,p)\\
 &=&2\langle\overrightarrow{y_mTx_m},\overrightarrow{y_mp}\rangle.
\end{eqnarray}
Also, it follows from Lemma \ref{coffi quasi lem} that
$ \langle\overrightarrow{y_mTx_m},\overrightarrow{y_mp}\rangle \leqslant \alpha_m \langle\overrightarrow{u_mTx_m},\overrightarrow{y_mp}\rangle$. Hence, we have
\begin{eqnarray}
\nonumber d^2(x_m,p)&=&\langle\overrightarrow{x_mp},\overrightarrow{x_mp}\rangle= \langle\overrightarrow{x_mTx_m},\overrightarrow{x_mp}\rangle+\langle\overrightarrow{Tx_mp},\overrightarrow{x_mp}\rangle\\
\nonumber &\leqslant&  \langle\overrightarrow{y_mTx_m},\overrightarrow{y_mp}\rangle+\langle\overrightarrow{Tx_mp},\overrightarrow{x_mp}\rangle \leqslant  \alpha_m \langle\overrightarrow{u_mTx_m},\overrightarrow{y_mp}\rangle+\langle\overrightarrow{Tx_mp},\overrightarrow{x_mp}\rangle\\
\nonumber &=&  \alpha_m \langle\overrightarrow{u_mTx_m},\overrightarrow{y_mx_m}\rangle+ \alpha_m \langle\overrightarrow{u_mp},\overrightarrow{x_mp}\rangle+ (1-\alpha_m)\langle\overrightarrow{Tx_mp},\overrightarrow{x_mp}\rangle\\
\nonumber &\leqslant&  \alpha_m d (u_m, Tx_m) d( y_m, x_m)+ \alpha_m \langle\overrightarrow{u_mo},\overrightarrow{x_mp}\rangle+\alpha_m \langle\overrightarrow{op},\overrightarrow{x_mp}\rangle\\
\nonumber &&+ (1-\alpha_m) d( Tx_m, p) d( x_m, p)\\
\nonumber &\leqslant& 2\alpha_m M d( y_m, x_m)+ \alpha_m \|u_m\| d( x_m, p) +\alpha_m \langle\overrightarrow{op},\overrightarrow{x_mp}\rangle+ (1-\alpha_m) d^2( x_m, p),
\end{eqnarray}
which implies that
\begin{eqnarray}\label{strong by weak}
 d^2(x_m,p)\leqslant  2 M d( y_m, x_m)+ \|u_m\| ( \|p\|+ M) +\langle\overrightarrow{op},\overrightarrow{x_mp}\rangle.
\end{eqnarray}
Since $\{x_m\}$ is bounded, by Lemma \ref{delta convergent subseq}, there exists a subsequence $\{x_{m_i}\}$ of $\{x_m\}$ which $\Delta$-converges to a point $q$. By Lemma \ref{deta limit is in C}, $q\in C$. It follows from (\ref{d(x_m, Tx_m) to 0}) and Lemma \ref{demiclosed of T} that $q\in F(T)$. Substituting $m_i$ and $q$, respectively, for $m$ and $p$ in (\ref{strong by weak}), we get
\begin{eqnarray}\label{strong by weak2}
 d^2(x_{m_i}, q)\leqslant  2 M d( y_{m_i}, x_{m_i})+ \|u_{m_i}\| ( \|q\|+ M) +\langle\overrightarrow{oq},\overrightarrow{x_{m_i}q}\rangle.
\end{eqnarray}
This together with Lemma \ref{delta iff limsup}, (\ref{d(x_m, y_m) to 0}) and $\Delta$-convergence of $\{x_{m_i}\}$ to $q$ implies that $x_{m_i}\to  q$ strongly as $i\to \infty$.\\
%That is every subsequence of $\{x_m\}$ has a convergent subsequence.\\
Now, if $\{x_{m_j}\}$ is a subsequence of $\{x_m\}$ which converges to a point $q'\in C$, then by using the same argument as in proof above, we get $q'\in F(T)$. Utilizing (\ref{strong by weak}), we have
\begin{eqnarray}
\nonumber  d^2(x_{m_i},q')\leqslant  2 M d( y_{m_i}, x_{m_i})+\|u_{m_i}\| ( \|q'\|+ M) +\langle\overrightarrow{oq'},\overrightarrow{x_{m_i}q'}\rangle
\end{eqnarray}
and
\begin{eqnarray}
\nonumber d^2(x_{m_j},q)\leqslant 2 M d( y_{m_j}, x_{m_j})+  \|u_{m_j}\| ( \|q\|+ M) +\langle\overrightarrow{oq},\overrightarrow{x_{m_j}q}\rangle.
\end{eqnarray}
By (\ref{d(x_m, y_m) to 0}) and continuity of $d(\cdot, x)$ for all $x\in X$, we get
\begin{eqnarray}\label{nearest}
d^2(q,q')\leqslant \langle\overrightarrow{oq'},\overrightarrow{qq'}\rangle\ \ \mbox{ and}\ \ \ d^2(q',q)\leqslant \langle\overrightarrow{oq},\overrightarrow{q'q}\rangle.
\end{eqnarray}
 Therefore, we obtain
\begin{eqnarray}
\nonumber 2d^2(q,q') \leqslant  \langle\overrightarrow{oq'},\overrightarrow{qq'}\rangle + \langle\overrightarrow{oq},\overrightarrow{q'q}\rangle
=\langle\overrightarrow{q'o},\overrightarrow{q'q}\rangle + \langle\overrightarrow{oq},\overrightarrow{q'q}\rangle=\langle\overrightarrow{q'q},\overrightarrow{q'q}\rangle=d^2(q,q').
\end{eqnarray}
Thus, $q'=q$. This shows that $\{x_m\}$ converges to $q\in F(T)$ and the proof is completed. $\ \ \ \Box$
%%%%%%%%%%%%%%%%%%%%%%%%%%%%%%%%%%%%%%%%%%%%%%%%%%%%%%%%%%%%%%%%%%%%%%%%%%%%%%%%%%%%%%%%%%%%%%%%%%%%%%%%%%%%%%%%%%%%%%%%%%%%%%%%%%%

\subsection{Convergence of an explicit algorithm\\}

\par
In this subsection we study strong convergence of an explicit algorithm to a fixed point of nonexpansive mappings.
\begin{algorithm}
Let $C$ be a nonempty closed convex subset of a Hadamard space $X$. Let $T : C \rightarrow C$ be a nonexpansive mapping. Define the iterative sequence $\{x_n\}$ as follows:
\begin{eqnarray}\label{explicit mann like}
\left\{\begin{array}{ll}
 x_0\in C,\ \mbox{chosen\ arbitrary},\\
         y_n=\alpha_n u_{n}\oplus(1-\alpha_n)Tx_n,  &  n\geqslant 0 , \\
         x_{n+1}=(1-\beta_n)x_n\oplus \beta_nP_Cy_n, &
       \end{array}
\right.
\end{eqnarray}
where $\{\alpha_n\}$  and $\{\beta_n\}$ are two sequences in $(0, 1)$, and the sequence $\{u_n\}\subset X$ is a perturbation
for the $n$-step iteration.

\end{algorithm}

%%%%%%%%%%%%%%%%%%%%%%%%%%%%%%%%%%%%%%%%%%%%%%%%%%%%%%%%%%%%%%%%%%%%%%%%%%%%%%%%%%%%%%%%%%%%%%%%%%%%%%%%%%%%%%%%%%%%%%%%%%%%%%%%%%%
\begin{theorem}
Let $F(T)\not=\emptyset$.
If the conditions (i) $\lim_{n\to \infty} \alpha_n=0$ and
$\sum_{n=0}^\infty \alpha_n=\infty$, (ii) $0< \liminf_{n\to
\infty} \beta_n\leqslant \limsup_{n\to \infty} \beta_n< 1$, (iii)
$\sum_{n=0}^\infty \alpha_n\|u_n\|<\infty$ satisfy, then the
sequence $\{x_n\}$ generated by the explicit method (\ref{explicit
mann like}) converges to a $q\in F(T)$.
\end{theorem}
\textbf{Proof.}
Let $p\in F(T)$. By Lemma \ref{triangle ine} and nonexpansiveness of $P_C$, we have
\begin{eqnarray}
\nonumber d(x_{n+1},p)&=&d\left((1-\beta_n)x_n\oplus \beta_nP_Cy_n,p\right)\\
\nonumber&\leqslant& (1-\beta_n)d(x_n,p ) + \beta_n d(P_Cy_n, p)\\
\nonumber&\leqslant& (1-\beta_n)d(x_n,p ) + \beta_n [d(\alpha_n u_{n}\oplus(1-\alpha_n)Tx_n , p)]\\
\nonumber&\leqslant& (1-\beta_n)d(x_n,p ) + \beta_n [\alpha_n d(u_{n}, p)+ (1-\alpha_n) d(Tx_n , p)]\\
\nonumber&\leqslant& (1-\beta_n)d(x_n,p ) + \beta_n [\alpha_n (\|u_n\|+\|p\|)+ (1-\alpha_n) d(x_n , p)]\\
\nonumber&\leqslant& (1-\alpha_n\beta_n)d(x_n,p ) + \beta_n\alpha_n \|p\|+ \alpha_n\|u_n\|\\
\nonumber&\leqslant& \max\{d(x_n,p ), \|p\|\} + \alpha_n\|u_n\|.
\end{eqnarray}
By induction, we get
\begin{eqnarray}
\nonumber d(x_{n+1},p) \leqslant \max\{d(x_0,p ), \|p\|\} + \sum_{i=0}^n\alpha_i\|u_i\|,
\end{eqnarray}
which together with condition (iii) implies that $\{x_n\}$ is bounded, so is the sequence $\{Tx_n\}\subset C$. Next, we prove that
\begin{eqnarray}\label{d(x_n+1 , x_n)}
 \lim_{n\to \infty}d(x_{n+1},x_n)=0.
\end{eqnarray}
Let $z_n=P_Cy_n$ for all $n\geqslant0$. It follows from nonexpansiveness of $P_C$ and Lemma \ref{hyberbolic ineqal} that
\begin{eqnarray}
\nonumber d(z_{n+1}, z_n)&\leqslant&d(y_{n+1}, y_n)=d\left(\alpha_{n+1} u_{n+1}\oplus(1-\alpha_{n+1})Tx_{n+1} , \alpha_n u_{n}\oplus(1-\alpha_n)Tx_n\right)\\
\nonumber&\leqslant& d\left(\alpha_{n+1} u_{n+1}\oplus(1-\alpha_{n+1})Tx_{n+1} , \alpha_{n+1} u_{n+1}\oplus(1-\alpha_{n+1})Tx_{n} \right)\\
\nonumber && + d \left( \alpha_{n+1} u_{n+1}\oplus(1-\alpha_{n+1})Tx_{n}, Tx_n \right) + d \left( Tx_n,  \alpha_n u_{n}\oplus(1-\alpha_n)Tx_n\right)\\
\nonumber&\leqslant& (1-\alpha_{n+1}) d (Tx_{n+1}, Tx_n) + \alpha_{n+1} d (u_{n+1}, Tx_n) + \alpha_{n} d (u_{n}, Tx_n)\\
\nonumber&\leqslant& (1-\alpha_{n+1}) d (x_{n+1}, x_n) + \alpha_{n+1} \|u_{n+1}\| + \alpha_{n} \|u_{n}\| + (\alpha_{n+1}+ \alpha_{n})\|Tx_n\|.
\end{eqnarray}
Hence,
\begin{eqnarray}
\nonumber d(z_{n+1}, z_n)- d (x_{n+1}, x_n ) \leqslant  \alpha_{n+1} \|u_{n+1}\|+ \alpha_{n} \|u_{n}\| + (\alpha_{n+1}+ \alpha_{n})\|Tx_n\|.
\end{eqnarray}
This together with (i) and (iii) implies that $\limsup_{n\to \infty} (d(z_{n+1}, z_n)- d (x_{n+1}, x_n ) )\leqslant 0$. It follows from Lemma \ref{suzuki lem} that
\begin{eqnarray}\label{d(z_n, x_n) to 0}
 \lim_{n\to \infty} d( z_n, x_n)=0.
\end{eqnarray}
Since $ d(x_{n+1}, x_n)= \beta_{n} d\left( z_n ,x_n\right) $, we get (\ref{d(x_n+1 , x_n)}). Now, we show that
\begin{eqnarray}\label{d(x_n , Tx_n)}
 \lim_{n\to \infty} d(x_n, Tx_n)=0.
\end{eqnarray}
Utilizing (\ref{explicit mann like}) and Lemma \ref{triangle ine}, we have
\begin{eqnarray}
\nonumber d(x_n, Tx_n)&\leqslant& d(x_n, x_{n+1}) + d( x_{n+1}, Tx_n)\\
\nonumber&\leqslant& d(x_n, x_{n+1}) + (1- \beta_n)d( x_{n}, Tx_n)+ \beta_n d(P_Cy_n, Tx_n)\\
\nonumber&\leqslant& d(x_n, x_{n+1}) + (1- \beta_n)d( x_{n}, Tx_n)+ \beta_n d(y_n, Tx_n)\\
\nonumber&=& d(x_n, x_{n+1}) + (1- \beta_n)d( x_{n}, Tx_n)+ \beta_n \alpha_nd(u_n, Tx_n).
\end{eqnarray}
Hence,
\begin{eqnarray}
\nonumber d(x_n, Tx_n)&\leqslant& \frac{1}{\beta_n} d(x_n, x_{n+1}) +\alpha_n (\|u_n\|+ \|Tx_n\|).
\end{eqnarray}
This together with (\ref{d(x_n+1 , x_n)}) and conditions (i)-(iii) implies that $\lim_{n\to \infty} d(x_n, Tx_n)=0$.  Moreover,
\begin{eqnarray}\label{d(x_n, y_n) to 0}
\nonumber d( x_n, y_n) &\leqslant& \alpha_n d( x_n, u_n)+ (1-\alpha_n)d( x_n, Tx_n)\\
\nonumber  &\leqslant& \alpha_n \|x_n\|+ \alpha_n \|u_n\| + (1-\alpha_n)d( x_n, Tx_n)\to 0\ \ \ (as\ n\to \infty).
\end{eqnarray}
Therefore, the inequality $d(z_n,y_n)\leqslant d(z_n,x_n)+d(x_n, y_n)$ together with (\ref{d(z_n, x_n) to 0}) implies that
\begin{eqnarray}\label{d(z_n, y_n) to 0}
\lim_{n\to \infty} d( z_n, y_n)=0.
\end{eqnarray}
\par
Let $\{x'_m\}$ be the sequence defined by the implicit method (\ref{mann like}), $q= \lim_{m\to \infty} x'_m$ and $y'_m=\alpha'_m u_{m}\oplus(1-\alpha'_m)Tx'_m$ for all $m\geqslant1$ where $\{\alpha'_m\}\subseteq (0,1)$ and $\{u'_m\}\subseteq X$ with $\lim_{m\to \infty} \alpha'_m=\lim_{m\to \infty} \|u'_m\|=0$. We show that $\limsup_{n\to \infty} \langle\overrightarrow{qo},\overrightarrow{qx_n}\rangle\leqslant0$. Similar to (\ref{x_m to y_m}) we may obtain that  $\langle\overrightarrow{x'_mTx'_m},\overrightarrow{x'_mx_n}\rangle\leqslant\langle\overrightarrow{y'_mTx'_m},\overrightarrow{y'_mx_n}\rangle$. Also, by Lemma \ref{coffi quasi lem}

\begin{eqnarray}
\nonumber d^2(x'_m,x_n)&=&\langle\overrightarrow{x'_mx_n},\overrightarrow{x'_mx_n}\rangle= \langle\overrightarrow{x'_mTx'_m},\overrightarrow{x'_mx_n}\rangle+\langle\overrightarrow{Tx'_mx_n},\overrightarrow{x'_mx_n}\rangle\\
\nonumber &\leqslant&  \langle\overrightarrow{y'_mTx'_m},\overrightarrow{y'_mx_n}\rangle+\langle\overrightarrow{Tx'_mx_n},\overrightarrow{x'_mx_n}\rangle\\
\nonumber &\leqslant&  \alpha'_m \langle\overrightarrow{u'_mTx'_m},\overrightarrow{y'_mx_n}\rangle+\langle\overrightarrow{Tx'_mx_n},\overrightarrow{x'_mx_n}\rangle\\
\nonumber &=&  \alpha'_m \langle\overrightarrow{u'_mTx'_m},\overrightarrow{y'_mx'_m}\rangle+ \alpha'_m \langle\overrightarrow{u'_mx_n},\overrightarrow{x'_mx_n}\rangle+ (1-\alpha'_m)\langle\overrightarrow{Tx'_mx_n},\overrightarrow{x'_mx_n}\rangle\\
\nonumber &=&  \alpha'_m \langle\overrightarrow{u'_mTx'_m},\overrightarrow{y'_mx'_m}\rangle+ \alpha'_m \langle\overrightarrow{u'_mx'_m},\overrightarrow{x'_mx_n}\rangle + \alpha'_m \langle\overrightarrow{x'_mx_n},\overrightarrow{x'_mx_n}\rangle\\
\nonumber &&+ (1-\alpha'_m)\langle\overrightarrow{Tx'_mTx_n},\overrightarrow{x'_mx_n}\rangle+(1-\alpha'_m)\langle\overrightarrow{Tx_nx_n},\overrightarrow{x'_mx_n}\rangle\\
\nonumber &\leqslant&  \alpha'_m d (u'_m, Tx'_m) d( y'_m, x'_m)+ \alpha'_m \langle\overrightarrow{u'_mo},\overrightarrow{x'_mx_n}\rangle+\alpha'_m \langle\overrightarrow{ox'_m},\overrightarrow{x'_mx_n}\rangle\\
\nonumber && + \alpha'_m d^2 (x'_m , x_n) + (1-\alpha'_m) d( Tx'_m, Tx_n) d( x'_m, x_n)\\
\nonumber && + (1-\alpha'_m) d( Tx_n, x_n) d( x'_m, x_n)\\
\nonumber &\leqslant& 2\alpha'_m M d( y'_m, x'_m)+ \alpha'_m \|u'_m\| d( x'_m, x_n) +\alpha'_m \langle\overrightarrow{ox'_m},\overrightarrow{x'_mx_n}\rangle\\
\nonumber &&+ d^2( x'_m, x_n)+ d( Tx_n, x_n) d( x'_m, x_n),
\end{eqnarray}
where $M>0$ is such that $\|u'_m\|, \|x'_m\|, \|Tx'_m\|\leqslant M$ for all $m\geqslant 1$. It follows that
\begin{eqnarray}
\nonumber \langle\overrightarrow{x'_mo},\overrightarrow{x'_mx_n}\rangle \leqslant 2M d( y'_m, x'_m) + \|u'_m\|M'+\frac{d( Tx_n, x_n)M'}{\alpha'_m},
\end{eqnarray}
where $M'>0 $ such that $d(x'_m, x_n)\leqslant M'$ for all $m\geqslant1$ and $n\geqslant 0$. It follows from (\ref{d(x_m, y_m) to 0}), $\|u'_m\|\to 0$ and (\ref{d(x_n , Tx_n)}) that
\begin{eqnarray}\label{lim m n}
 \limsup_{m\to \infty} \limsup_{n\to \infty}\langle\overrightarrow{x'_mo},\overrightarrow{x'_mx_n}\rangle \leqslant 0.
\end{eqnarray}
We note that
\begin{eqnarray}
\nonumber \langle\overrightarrow{qo},\overrightarrow{qx_n}\rangle &=& \langle\overrightarrow{qo},\overrightarrow{qx'_m}\rangle+ \langle\overrightarrow{qx'_m},\overrightarrow{x'_mx_n}\rangle + \langle\overrightarrow{x'_mo},\overrightarrow{x'_mx_n}\rangle\\
\nonumber &\leqslant&  \langle\overrightarrow{qo},\overrightarrow{qx'_m}\rangle+ d (q, x'_m) d(x'_m, x_n) + \langle\overrightarrow{x'_mo},\overrightarrow{x'_mx_n}\rangle\\
\nonumber &\leqslant&  \langle\overrightarrow{qo},\overrightarrow{qx'_m}\rangle+ d (q, x'_m)M' + \langle\overrightarrow{x'_mo},\overrightarrow{x'_mx_n}\rangle.
\end{eqnarray}
This together with (\ref{lim m n}) and $\lim_{m\to \infty} x'_m=q$ implies that
\begin{eqnarray}\label{lim m n qx_n}
  \limsup_{n\to \infty} \langle\overrightarrow{qo},\overrightarrow{qx_n}\rangle= \limsup_{m\to \infty}\limsup_{n\to \infty} \langle\overrightarrow{qo},\overrightarrow{qx_n}\rangle \leqslant \limsup_{m\to \infty} \limsup_{n\to \infty}\langle\overrightarrow{x'_mo},\overrightarrow{x'_mx_n}\rangle \leqslant 0.
\end{eqnarray}
Since
\begin{eqnarray}
\nonumber \langle\overrightarrow{qo},\overrightarrow{qz_n}\rangle = \langle\overrightarrow{qo},\overrightarrow{qx_n}\rangle+ \langle\overrightarrow{qo},\overrightarrow{x_nz_n}\rangle \leqslant  \langle\overrightarrow{qo},\overrightarrow{qx_n}\rangle+ \|q\| d(x_n, z_n)
\end{eqnarray}
and
\begin{eqnarray}
\nonumber \langle\overrightarrow{qo},\overrightarrow{qTx_n}\rangle = \langle\overrightarrow{qo},\overrightarrow{qx_n}\rangle+ \langle\overrightarrow{qo},\overrightarrow{x_nTx_n}\rangle \leqslant  \langle\overrightarrow{qo},\overrightarrow{qx_n}\rangle+ \|q\| d(x_n, Tx_n),
\end{eqnarray}
using (\ref{d(z_n, x_n) to 0}), (\ref{d(x_n , Tx_n)}) and (\ref{lim m n qx_n}), we have
\begin{eqnarray}\label{lim m n qx_n, qTx_n}
  \limsup_{n\to \infty} \langle\overrightarrow{qo},\overrightarrow{qz_n}\rangle\leqslant 0\ \ \ \ \mbox{and}\ \ \ \ \limsup_{n\to \infty} \langle\overrightarrow{qo},\overrightarrow{qTx_n}\rangle\leqslant 0.
\end{eqnarray}
Finally, we show that $\lim_{n\to \infty}x_n =q$. Since $z_n=P_Cy_n$ and $q\in C$, it follows from Theorem \ref{proj} that   $\langle\overrightarrow{z_ny_n},\overrightarrow{z_nq}\rangle\leqslant 0$. Using (\ref{explicit mann like}) and Lemma \ref{para ine} and \ref{coffi quasi lem}, we have
\begin{eqnarray}\label{x_n+1, q}
\nonumber d^2(x_{n+1},q)&\leqslant& (1-\beta_n)d^2(x_{n},q) + \beta_n d^2(z_{n},q)\\
\nonumber &=& (1-\beta_n)d^2(x_{n},q) + \beta_n \langle\overrightarrow{z_ny_n},\overrightarrow{z_nq}\rangle
+  \beta_n\langle\overrightarrow{y_nq},\overrightarrow{z_nq}\rangle\\
\nonumber &\leqslant& (1-\beta_n)d^2(x_{n},q) +  \beta_n\langle\overrightarrow{y_nq},\overrightarrow{z_nq}\rangle\\
\nonumber &=& (1-\beta_n)d^2(x_{n},q) +  \beta_n\langle\overrightarrow{y_nTx_n},\overrightarrow{z_nq}\rangle + \beta_n\langle\overrightarrow{Tx_nq},\overrightarrow{z_nq}\rangle\\
\nonumber &=& (1-\beta_n)d^2(x_{n},q) +  \beta_n\langle\overrightarrow{y_nTx_n},\overrightarrow{z_ny_n}\rangle + \beta_n\langle\overrightarrow{y_nTx_n},\overrightarrow{y_nq}\rangle \\
\nonumber && + \beta_n\langle\overrightarrow{Tx_nq},\overrightarrow{z_nq}\rangle\\
\nonumber &\leqslant& (1-\beta_n)d^2(x_{n},q) +  \beta_n\langle\overrightarrow{y_nTx_n},\overrightarrow{z_ny_n}\rangle + \beta_n\alpha_n \langle\overrightarrow{u_nTx_n},\overrightarrow{y_nq}\rangle \\
\nonumber &&+ \beta_n\langle\overrightarrow{Tx_nq},\overrightarrow{z_nq}\rangle\\
\nonumber &=& (1-\beta_n)d^2(x_{n},q) +  \beta_n\langle\overrightarrow{y_nTx_n},\overrightarrow{z_ny_n}\rangle
+ \beta_n\alpha_n \langle\overrightarrow{u_nTx_n},\overrightarrow{y_nz_n}\rangle \\
\nonumber &&+ \beta_n \left(\alpha_n \langle\overrightarrow{u_nq},\overrightarrow{z_nq}\rangle +  (1-\alpha_n)\langle\overrightarrow{Tx_nq},\overrightarrow{z_nq}\rangle\right)\\
\nonumber &\leqslant & (1-\beta_n)d^2(x_{n},q) +  \beta_n d (y_n, Tx_n) d (z_n, y_n) + \beta_n\alpha_n d(u_n, Tx_n) d (z_n, y_n) \\
\nonumber &&+ \beta_n \left(\alpha_n \langle\overrightarrow{u_no},\overrightarrow{z_nq}\rangle + \alpha_n \langle\overrightarrow{oq},\overrightarrow{z_nq}\rangle +  (1-\alpha_n)\langle\overrightarrow{Tx_nq},\overrightarrow{z_nq}\rangle\right)\\
\nonumber &= & (1-\beta_n)d^2(x_{n},q) + 2\beta_n\alpha_n d(u_n, Tx_n) d (z_n, y_n) \\
 &&+ \beta_n \left(\alpha_n \langle\overrightarrow{u_no},\overrightarrow{z_nq}\rangle + \alpha_n \langle\overrightarrow{oq},\overrightarrow{z_nq}\rangle +  (1-\alpha_n)\langle\overrightarrow{Tx_nq},\overrightarrow{z_nq}\rangle\right).
\end{eqnarray}
We note that
\begin{eqnarray}\label{Tx_nq, z_nq}
\langle\overrightarrow{Tx_nq},\overrightarrow{z_nq}\rangle\leqslant d(Tx_n, q) d(z_n, q)\leqslant d(x_n, q) d(y_n, q)\leqslant \frac{1}{2}\left(d^2(x_n, q)+ d^2(y_n, q)\right).
\end{eqnarray}
Also, using Lemma \ref{hyberbolic ineqal} and \ref{coffi quasi lem 2} we have
\begin{eqnarray}\label{d^2(y_n,q)}
\nonumber d^2(y_n, q) &\leqslant& \left( d(y_n,  \alpha_n o \oplus (1-\alpha_n)Tx_n ) + d(\alpha_n o \oplus (1-\alpha_n)Tx_n , q)\right)^2\\
\nonumber &\leqslant& \left( \alpha_n \|u_n\| + d(\alpha_n o \oplus (1-\alpha_n)Tx_n , q)\right)^2 \\
\nonumber &=& \alpha_n^2 \|u_n\|^2 + d^2(\alpha_n o \oplus (1-\alpha_n)Tx_n , q) + 2 \alpha_n \|u_n\|d(\alpha_n o \oplus (1-\alpha_n)Tx_n , q)\\
\nonumber &\leqslant&  \alpha_n^2 \|q\|^2 + (1-\alpha_n)^2 d^2(Tx_n , q)+ 2 \alpha_n (1- \alpha_n) \langle\overrightarrow{oq},\overrightarrow{Tx_nq}\rangle \\
\nonumber &&+ \alpha_n^2 \|u_n\|^2 + 2 \alpha_n \|u_n\|d(\alpha_n o \oplus (1-\alpha_n)Tx_n , q)\\
\nonumber &\leqslant&  (1-\alpha_n) d^2(x_n , q)+  \alpha_n \left(\alpha_n \|q\|^2 + 2  (1- \alpha_n) \langle\overrightarrow{oq},\overrightarrow{Tx_nq}\rangle\right)\\
 &&+ \alpha_n \|u_n\| M'',
\end{eqnarray}
where $M''=\sup_n\{ \alpha_n \|u_n\| + 2 d((1-\alpha_n)Tx_n \oplus \alpha_n o , q)\}$. Using (\ref{x_n+1, q}), (\ref{Tx_nq, z_nq}) and (\ref{d^2(y_n,q)}), we obtain

\begin{eqnarray}
\nonumber d^2(x_{n+1},q)&\leqslant& (1-\beta_n)d^2(x_{n},q) + 2\beta_n\alpha_n d(u_n, Tx_n) d (z_n, y_n) \\
\nonumber &&+ \beta_n \left(\alpha_n \langle\overrightarrow{u_no},\overrightarrow{z_nq}\rangle + \alpha_n \langle\overrightarrow{oq},\overrightarrow{z_nq}\rangle +  \frac{(1-\alpha_n)}{2}\left(d^2(x_n, q)+ d^2(y_n, q)\right)\right)\\
\nonumber&\leqslant&  (1-\beta_n)d^2(x_{n},q) + 2\beta_n\alpha_n d(u_n, Tx_n) d (z_n, y_n) \\
\nonumber &&+ \beta_n \alpha_n \|u_n\| d(z_n,q) + \beta_n \alpha_n \langle\overrightarrow{oq},\overrightarrow{z_nq}\rangle \\
\nonumber &&+  \frac{(1-\alpha_n)\beta_n}{2} \left(d^2(x_n, q)+ (1-\alpha_n) d^2(x_n , q)\right)\\
\nonumber &&+  \frac{(1-\alpha_n)\beta_n}{2} \left( \alpha_n \left(\alpha_n \|q\|^2 + 2  (1- \alpha_n) \langle\overrightarrow{oq},\overrightarrow{Tx_nq}\rangle\right)+ \alpha_n \|u_n\| M''\right)\\
\nonumber&\leqslant&  (1-\beta_n \alpha_n)d^2(x_{n},q) + \beta_n\alpha_n \left( 2 \|Tx_n\| d (z_n, y_n) + \langle\overrightarrow{oq},\overrightarrow{z_nq}\rangle + \alpha_n \|q\|^2\right. \\
\nonumber &&\left. + 2  (1- \alpha_n) \langle\overrightarrow{oq},\overrightarrow{Tx_nq}\rangle\right)+ \alpha_n \|u_n\| \left(2 d (z_n, y_n) +d(z_n,q) + M''\right)\\
\nonumber&=& (1-\gamma_n )d^2(x_{n},q) + \gamma_n\delta_n + \sigma_n,
\end{eqnarray}
where
$\gamma_n= \beta_n \alpha_n$, $\delta_n= 2 \|Tx_n\| d (z_n, y_n) + \langle\overrightarrow{oq},\overrightarrow{z_nq}\rangle +\alpha_n \|q\|^2 + 2  (1- \alpha_n) \langle\overrightarrow{oq},\overrightarrow{Tx_nq}\rangle $ and $\sigma_n= \alpha_n \|u_n\| (2 d (z_n, y_n)+d(z_n,q) + M'')$. Now, considering conditions (i)-(iii), (\ref{d(z_n, y_n) to 0}), (\ref{lim m n qx_n, qTx_n}) and applying Lemma \ref{Liu' lem} to the last inequality, we conclude that
$\lim_{n\to \infty}x_n=q$. This completes the proof.
$\ \ \ \ \Box$
\begin{remark}
The algorithms (\ref{mann like}) and (\ref{explicit mann like}) converge strongly to $P_{F(T)}o$, the nearest point of $F(T)$ to $o$. In fact, a similar method as in proof of (\ref{nearest}) shows that $d^2(q,p)\leqslant \langle\overrightarrow{op},\overrightarrow{qp}\rangle$ for all $p\in F(T)$. Which is equivalent to
$ \langle\overrightarrow{qo},\overrightarrow{qp}\rangle \leqslant 0$ for all $p\in F(T)$. It follows from Theorem \ref{proj} that $q=P_{F(T)}o$.
\end{remark}
\textbf{Acknowledgement.}
The authors are grateful to Professor N. Shahzad for his valuable suggestions and the anonymous referee for his/her careful reading of the paper.\\\\
%% appendix sections are then done as normal sections
%% \appendix

%% \section{}
%% \label{}

%% References
%%
%% Following citation commands can be used in the body text:
%% Usage of \cite is as follows:
%% \cite{key}   ==>> [#]
%% \cite[chap. 2]{key} ==>> [#, chap. 2]
%% \citet{key}   ==>> Author [#]
%% References with bibTeX database:
%\bibliographystyle{model3a-num-names}
%\bibliography{<your-bib-database>}
%\noindent\textbf{References}

%% Authors are advised to submit their bibtex database files. They are
%% requested to list a bibtex style file in the manuscript if they do
%% not want to use model3a-num-names.bst.

%% References without bibTeX database:

% \begin{thebibliography}{00}

%% \bibitem must have the following form:
%% \bibitem{key}\ldots
%%

% \bibitem{}

% \end{thebibliography}

\end{document}